# ANALYTIC CENTER CUTTING PLANE METHODS FOR VARIATIONAL INEQUALITIES OVER CONVEX BODIES

**Renying Zeng**[1]

School of Mathematical Sciences, Chongqing Normal University, Chongqing, China

**ABSTRACT**

An analytic center cutting plane method is an iterative algorithm based on the computation of analytic centers. In this paper, we propose some analytic center cutting plane methods for solving quasimonotone or pseudomonotone variational inequalities whose domains are bounded or unbounded convex bodies.

## 1. Introduction and Preliminaries

Some recent developments in solving variational inequalities are analytic center cutting plane methods. An analytic center cutting plane method is an interior algorithm based on the computation of analytic centers. In order to work with analytic center cutting plane methods, some authors assume that the feasible sets of variational inequalities are polytopes, e.g., see [2, 3, 5, 6], while others pay more attention to problems with infinitely many linear constraints, e.g., see [4, 10], etc. Analytic center cutting plane methods also can be used to other types of optimization problems, like mathematical programming with equilibrium constraints [9], convex programming [16, 23], conic programming [17], stochastic programming [18, 19], and combinatorial optimization [23]. In this paper, we propose some analytic center cutting plane methods for solving pseudomonotone or quasimonotone variational inequalities whose feasible sets are bounded or unbounded convex bodies (including the n-dimensional Euclidean space $\mathbb{R}^n$ itself).

---

[1] Mailing address: Mathematics Department, Saskatchewan Polytechnic, 1130 Idylwyld Dr. N, Saskatoon, Saskatchewan, Canada S7L 4J7. Email: renying.zeng@saskpolytech.ca

Let $X$ be a non-empty subset of the n-dimensional Euclidean space $\mathbb{R}^n$ and let $F: X \to \mathbb{R}^n$ be a function. We call that a point $x^* \in X$ is a solution of the **variational inequality** $VI(F, X)$ if

$$F(x^*)^T (x - x^*) \geq 0, \forall x \in X. \tag{1}$$

The point $x^* \in X$ is a solution of the **dual variational inequality** $VID(F, X)$ if

$$F(x)^T (x - x^*) \geq 0, \forall x \in X. \tag{2}$$

We denote by $X^*$ the set of solutions of $VI(F, X)$, and by $X_D^*$ the set of solutions of $VID(F, X)$.

Given $VI[F, X]$ ($VID[F, X]$), the gap function is defined as

$$g_X(x) = \max_{y \in X} F(x)^T (x - y), x \in X,$$

$$(f_X(x) = \max_{y \in X} F(y)^T (x - y), x \in X).$$

We note that $g_X(x) \geq 0$, $f_X(x) \geq 0$, and

$$\arg\min_{x \in X} g_X(x) = \{g_X(x) = 0; x \in X\}.$$

$$= \arg\min_{x \in X} f_X(x) = \{g_X(x) = 0; x \in X\}$$

Therefore, we have

**Lemma 2** A point $x^* \in X$ is a solution of $VI[F, X]$ ($VID[F, X]$) if and only if $g_X(x^*) = 0$ ($f_X(x^*) = 0$). □

A point $x^* \in X$ is said to be a **$\varepsilon$-solution** of the variational inequality (1) if $g_X(x^*) < \varepsilon$.

A function $F: X \to \mathbb{R}^n$ is said to be **monotone** on $X$ if

$$(F(y) - F(x))^T (y - x) \geq 0, \forall x, y \in X,$$

**strongly monotone** if

$$(F(y) - F(x))^T (y - x) \geq M \|y - x\|, \forall x, y \in X \text{ with a constant } M > 0,$$

**quasimonotone** on $X$ if

$$F(x)^T (y - x) > 0 \Rightarrow F(y)^T (y - x) \geq 0, \forall x, y \in X,$$

**pseudomonotone** on $X$ if

$$F(x)^T (y - x) \geq 0 \Rightarrow F(y)^T (y - x) \geq 0, \forall x, y \in X,$$

**strictly pseudomonotone** on $X$ if

$$F(x)^T (y - x) \geq 0 \Rightarrow F(y)^T (y - x) > 0, \forall x, y \in X,$$

**pseudomonotone plus** on $X$ if it is pseudomonotone on $X$ and if

$$\left. \begin{array}{l} F(x)^T (y - x) \geq 0 \\ F(y)^T (y - x) = 0 \end{array} \right\} \Rightarrow F(x) = F(y), \forall x, y \in X,$$

and, **strongly pseudomonotone** on $X$ if there exist constants $M > 0$, $\alpha > 0$, such that

$$F(x)^T (y - x) \geq 0 \Rightarrow F(y)^T (y - x) \geq M \|y - x\|^\alpha, \forall x, y \in X.$$

From Auslender [1] we have Lemma 1.

**Lemma 1** if $F$ is continuous, then a solution of $VID$ $(F, X)$ is a solution of $VI$ $(F, X)$; and if $F$ is continuous and strictly pseudomonotone, then $x^* \in X$ is a solution of $VI$ $(F, X)$ if and only if it is a solution of $VID$ $(F, X)$. □

## 2. Compact Convex Bodies

A **polytope** is a set $P \subseteq \mathbb{R}^n$ which is the convex hull of a finite set.

A **polyhedron** is a set

$$\{x \in \mathsf{R}^n; A^T x \le b\} \subseteq \mathsf{R}^n,$$

where $b \in \mathsf{R}^n$, and $A$ is an $m \times n$ matrix.

Every polytope is a polyhedron whereas not every polyhedron is a polytope.

H. Minkowski proved the following Lemma 3 in 1896.

**Lemma 3** A set $P \subseteq \mathsf{R}^n$ is a polytope if and only if it is a bounded polyhedron.□

We make the following assumptions for polytopes throughout this paper.

(**a**) **Interior Assumption**: A polytope is always a full-dimensional polytope and that includes $0 \le x \le e$, where $e$ is a vector of all ones.

We note that if a polytope has non-empty interior, then (**a**) can be met by re-scaling.

A **convex body** $X \subseteq \mathsf{R}^n$ is a convex and bounded subset with non-empty interior.

A **rectangle** $B \subseteq \mathsf{R}^n$ is defined by

$$B = \{x = (x_1, x_2, \cdots, x_n) \in \mathsf{R}^n; a_i \le x_i \le b_i\}$$

where $a_i, b_i \in \mathsf{R}$.

A rectangle can also be given by some inequalities

$$B = \{x \in \mathsf{R}^n; H^T x \le b\},$$

where $H^T x = b$ is a finite set of hyperplanes, $H$ is an $m \times n$ matrix. And, if we denote by $V$ the finite set of all vertices of $B$, then

$$B = con\,(V).$$

**Theorem 1** A bounded subset $X \subseteq \mathbb{R}^n$ is a compact convex body if and only if there exists a sequence of polytopes $\{C_j\}$, which satisfies $C_j \subseteq C_{j+1}$ ($j = 1, \cdots$), such that
$$(\bigcup_{j=1}^{\infty} C_j)^c = X.$$

**Proof.** The sufficiency is trivial. We only prove the necessity.

Since $X$ is bounded, there exists a rectangle $B$ such that $X \subseteq B$.

Take a partition $P_1$ of $B$. Then $B$ is divided into a set of finite subrectangles by a finite set of hyperplanes. Let $D_1 = \bigcup_{j=1}^{k_1} B_{1(j)}$, where $B_{1(j)}$ ($j = 1, \ldots, k_1$) are all the subrectangles which lie entirely within $X$. Let $V_1$ be the set of all vertices of $B_{1(j)}$ ($j = 1, \ldots, k_1$), then $V_1$ is a finite set. So, $C_1 = con(V_1)$ is a polytope, and it obviously satisfies
$$D_1 \subseteq C_1 \subseteq X.$$
(For the case of a 2-dimensional Euclidean space, see Figure 2.)

Take a finer partition $P_2$ of $B$. Similarly, we have a set $D_2 = \bigcup_{j=1}^{k_2} B_{2(j)}$, where $B_{2(j)}$ ($j = 1, \ldots, k_2$) are all the subrectangles which are corresponding to $P_2$ and lie entirely within $X$; and we have a polytope $C_2 = con(V_2)$, where $V_2$ is the set of all vertices of $B_{2(j)}$ ($j = 1, \ldots, k_2$), such that
$$C_1 \subseteq C_2 \subseteq X.$$

By mathematical induction, there exists a sequence of polytopes $\{C_j\}$ which satisfies
$$C_j \subseteq C_{j+1} \subseteq X \, (j = 1, \cdots).$$

It is easy to see that $(\bigcup_{j=1}^{\infty} C_j)^c = X$.

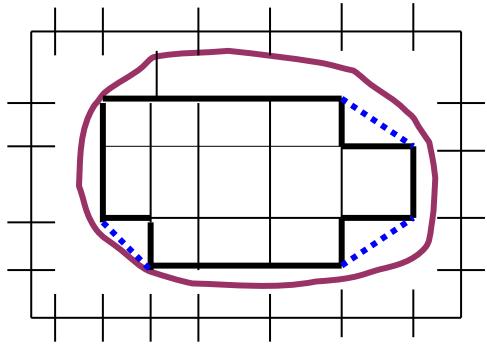

**Figure 2** $C_1 = con(V_1)$

It is quite straightforward to prove the following Corollary 1, Proposition 1, and Proposition 2.

**Corollary 1** A subset $X \subseteq \mathbb{R}^n$ is a compact convex body if there exists a is uniformly bounded sequence of polytopes $\{C_j\}$, i.e., $C_j \subseteq B$ for a given rectangle $B$, such that

$$(\bigcup_{j=1}^{\infty} C_j)^c = X. \square$$

**Proposition 1** Let $X \subseteq \mathbb{R}^n$ be a compact convex body and $F: X \to \mathbb{R}^n$ a continuous function, then the variational inequality $VI[F, X]$ has solutions. $\square$

**Proposition 2** Let $X \subseteq \mathbb{R}^n$ be a compact convex body and $F: X \to \mathbb{R}^n$ a continuous and strictly psuedomonotone function, then the variational inequality $VI[F, X]$ has a unique solution. $\square$

## 3. Generalized Analytic Center Cutting Plane Algorithms for Solving Pseudomonotone Variational Inequalities

For any polytope $\{x \in \mathbb{R}^n; A^T x \leq b\}$,

$$\{x \in \mathbb{R}^n; A^T x + s = b, s = (s_1, s_2, \cdots, s_n), s_i \geq 0\}$$

is associated with the potential function

$$\varphi = \sum_{i=1}^{n} \ln s_i.$$

It is known that an **analytic center** is the maximizer of the potential function $\varphi$, and the unique solution of the system

$$A^T y = 0$$
$$A^T x + s = b$$
$$Y^T s = e,$$

where $y$ is a positive dual vector, and $Y$ the diagonal matrix built upon $y$.

An **approximate analytic center** [12] is the maximizer of the potential function $\varphi$, and the unique solution of the system

$$A^T z = 0$$
$$A^T x + s = b$$
$$\| Z^T s - e \| \leq \eta < 1.$$

Where $z$ is a dual vector, and $Z$ the diagonal matrix built upon $z$.

Now we modify Goffin, Marcotte and Zhu's [3] Algorithm 1 to solve $VI(F, X)$. We propose an algorithm for solving variational inequalities, whose feasible sets are compact convex bodies.

From Theorem 1, we have a sequence of variational inequalities $VI[F, C_j]$ ($j = 1, \cdots$) induced by the original variational inequality $VI[F, X]$, where the polytope $C_j$ is given by the linear inequalities $A_j^T x = b_j$, $x, b_j \in \mathsf{R}^n$, and $A_j$ is an $m \times n$ matrix. So, we may apply the algorithm in [3] to each $VI[F, C_j]$. The following Algorithm 1 is using this idea, but the algorithm in [3] is applied to $VI[F, C_j]$ for only certain number of iterations until we get

$$g_{C_j}(x_j^k) < \frac{1}{2^j} \quad (j = 1, \cdots),$$

by applying Theorem 1 of [3].

**Algorithm 1.**

**Step 1.** (initialization)
$k = 0, j = 1, A^k = A_j, b^k = b_j, C_j^k = \{x \in \mathsf{R}^n ; A_j^k x \leq b_j^k\}$;

**Step 2**. (computation of an approximate analytic center)

Find an approximate analytic center $x_j^k$ of $C_j^k$;

**Step 3**. (stop criterion)

Compute $g_X(x_j^k) = \max_{x \in X} F(x_j^k)^T(x_j^k - x)$,

if $g_X(x_j^k) = 0$, **then** STOP,

        **else** GOTO step 4;

**Step 4**. (find an $\varepsilon$-solution for $\varepsilon = \dfrac{1}{2^j}$ )

Compute $g_{C_j}(x_j^k) = \max_{x \in C_j} F(x_j^k)^T(x_j^k - x)$,

if $g_{C_j}(x_j^k) < \dfrac{1}{2^j}$, **then** increase $j$ by 1 RETURN TO Step1,

        **else** GOTO Step 5;

**Step 5**. (cut generation)

Set $A_j^{k+1} = \begin{bmatrix} A_j^k \\ F(x_j^k)^T \end{bmatrix}$, $b_j^{k+1} = \begin{bmatrix} b_j^k \\ F(x_j^k)^T x_j^k \end{bmatrix}$,

$H_j^k = \{x \in \mathsf{R}^n; F(x_j^k)^T(x - x_j^k) = 0\}$ is the new cutting plane for $VI(F, C_j^k)$.

Increase $k$ by one GOTO Step 2.

**Theorem 2** Let $F: X \to \mathsf{R}^n$ be pseudomonotone plus on a compact convex body $X$, then Algorithm 1 either stops with a solution of $VI(F, X)$ after a finite number of iterations, or there exists a subsequence of the infinite sequence $\{x_j^k\}$ that converges to a point $x^* \in X^*$.

**Proof**. According to Algorithm 1 and Theorem 1 of [3], for any given $j$, after finite number of iterations, $\exists x_j \in C_j$ such that

$$g_{C_j}(x_j) < \frac{1}{2^j}.$$

Since $X$ is compact, there exists a subsequence $\{x_{j(q)}\}$ of $\{x_j\}$ and a point $x^* \in X$ such that

$$\lim_{q \to \infty} x_{j(q)} = x^*.$$

$\forall p < j$, we have

$$g_{C_p}(x_j) \le g_{C_j}(x_j) < \frac{1}{2^j}.$$

On the other hand, due to the compactness of $X$, $\exists N > 0$ such that $\|y\| \le N$, $\forall y \in X$. Since

$$\max_{y \in X} F(x')^T (x'-y) \ge 0, \text{ and } \max_{y \in X} F(x'')^T (x''-y) \ge 0,$$

for $\forall x', x'' \in X$,

$$|g_X(x') - g_X(x'')|$$
$$= |\max_{y \in X} F(x')^T (x'-y) - \max_{y \in X} F(x'')^T (x''-y)|$$
$$\le \max_{y \in X} |F(x')^T (x'-y) - F(x'')^T (x''-y)|$$
$$= \max_{y \in X} |[F(x')^T x' - F(x'')^T x''] + [F(x'')^T y - F(x')^T y]|$$
$$\le |[F(x')^T x' - F(x'')^T x'']| + \max_{y \in X} |F(x'')^T y - F(x')^T y|$$
$$\le |F(x')^T x' - F(x'')^T x''| + N \|F(x'') - F(x')\|.$$

By the continuities of $F(x)$ and $F(x)^T x$, $g_X(x)$ is a continuous function on $X$.

Consequently, $\forall p$

$$g_{C_p}(x^*) = \lim_{q \to \infty} g_{C_p}(x_{j(q)}) \le \lim_{q \to \infty} g_{C_{j(q)}}(x_{j(q)}) \le \lim_{q \to \infty} \frac{1}{2^{j(q)}} = 0.$$

Then, we have

$$g_{\bigcup_{j=1}^{\infty} C_j}(x^*) = \max_{y \in \bigcup_{j=1}^{\infty} C_j} F(x^*)(x^* - y)$$
$$\le \sum_{j=1}^{\infty} \max_{y \in C_j} F(x^*)(x^* - y) = \sum_{j=1}^{\infty} g_{C_j}(x^*) = 0.$$

On the other hand, $\forall y \in X$, $\exists \{y_i\} \subseteq \bigcup_{j=1}^{\infty} C_j$ such that

$$\lim_{i \to \infty} y_i = y.$$

Because

$$|F(x^*)(x^* - y_i)| \leq \max_{y \in \bigcup_{j=1}^{\infty} C_j} F(x^*)(x^* - y)$$

$$= g_{\bigcup_{j=1}^{\infty} C_j}(x^*) = 0,$$

we have

$$|F(x^*)(x^* - y)| = \lim_{i \to \infty} |F(x^*)(x^* - y_i)| = 0.$$

Therefore,

$$g_X(x^*) = \max_{y \in X} F(x^*)(x^* - y) = 0.$$

Which deduces that $x^*$ is a solution of $VI(F, X)$. □

Algorithm 1 usually generates an infinite sequence. In order to terminate at a finite number of iterations, we change the stop criterion, Step 3 in Algorithm 1, to get the following Algorithm 2.

**Algorithm 2.** Step 1, Step 2, Step 4, and Step 5 are same as those of Algorithm 1.

**Step 3**. (stop criterion)

Compute $g_X(x_j^k) = \max_{x \in X} F(x_j^k)^T(x_j^k - x)$,

if $g_X(x_j^k) < \varepsilon$, **then** STOP,

            **else** GOTO step 4.

From Theorem 2 we have

**Theorem 3** Let $F: X \to \mathbb{R}^n$ be pseudomonotone plus on a compact convex body $X$, then Algorithm 2 stops with a $\varepsilon$-solution of $VI(F, X)$ after a finite number of iterations. □

# 4. Generalized Analytic Center Cutting Plane Algorithms for Solving Quasimonotone Variational Inequalities

In this section, we are going to modify Marcotte and Zhu's [2] approach to solve quasimonotone variational inequalities $VI(F, X)$. We assume that the feasible sets are compact convex bodies.

From Theorem 1 there is a sequence of variational inequalities $VI[F, C_j]$ ($j = 1, \cdots$) induced by the original variational inequality $VI[F, X]$.

According to [2], the following are the elements required in the construction of algorithms for solving quasimonotone variational inequalities.

For any given $j$, let the auxiliary function $\Gamma_j(y,x): \mathbb{R}^n \to \mathbb{R}^n$ be continuous in $x$ and strong monotone in $y$, i.e.,

$$\langle \Gamma_j(y',x) - \Gamma_j(y'',x), y'-y'' \rangle \geq \beta_j \| y'-y'' \|^2, \forall y', y'' \in X,$$

for $\beta_j > 0$. $\beta_j$ is said to be the strong monotonicity constant for $\Gamma_j(y,x): \mathbb{R}^n \to \mathbb{R}^n$. The function $\Gamma_j$ is associated with the variational inequality $AVI[\Gamma, X, x]$ whose solution $w_j(x)$ satisfies

$$\langle \Gamma_j(w(x),x) - \Gamma_j(x,x) + F(x), y - w_j(x) \rangle \geq 0, \forall y \in X.$$

It is known that $w_j(x)$ are continuous [15], and that $x$ is a solution of $VI[F,C_j]$ if and only if it is a fixed point of $w$.

Let $\rho_j$ and $\alpha_j$ be positive numbers less than 1 and $\beta_j$, respectively. Let $l(j)$ (which depends on $x$) be the smallest nonnegative integer for which

$$\langle F(x + \rho_j^{l(j)}(w_j(x) - x)), x - w_j(x) \rangle \geq \alpha_j \| w_j(x) - x_j \|^2.$$

Define

$$G_j(x) = F(x + \rho_j^{l(j)}(w_j(x) - x)).$$

If $x_j^*$ is a solution of $VI[F,C_j]$, then $w_j(x_j^*) = x_j^*$, $l(j) = 0$, and $G_j(x_j^*) = F(x_j^*)$.

**Algorithm 3.**

**Step 1**. (initialization)

Let $\beta_j > 0$ be the strong monotonicity constant for $\Gamma_j(y,x): \mathsf{R}^n \to \mathsf{R}^n$, with respect to $y$ and let $\alpha_j \in (0, \beta_j)$.

$k = 0, j = 1, A^k = A_j, b^k = b_j, C_j^k = \{x \in \mathsf{R}^n; A_j^k x \leq b_j^k\}$;

**Step 2**. (computation of an approximate analytic center)

Find an approximate analytic center $x_j^k$ of $C_j^k$;

**Step 3**. (stop criterion)

Compute $g_X(x_j^k) = \max_{x \in X} F(x_j^k)^T (x_j^k - x)$,

if $g_X(x_j^k) = 0$, **then** STOP,

   **else** GOTO step 4;

**Step 4**. (find an $\varepsilon$-solution for $\varepsilon = \dfrac{1}{2^j}$)

Compute $g_{C_j}(x_j^k) = \max_{x \in C_j} F(x_j^k)^T (x_j^k - x)$,

if $g_{C_j}(x_j^k) < \dfrac{1}{2^j}$, **then** increase $j$ by 1 RETURN TO Step1,

   **else** GOTO Step 5;

**Step 5**. (auxiliary variational inequality)

Let $w_j(x_j^k)$ satisfies the variational inequality

$$\langle F(x_j^k) + \Gamma_j(w_j(x_j^k), x_j^k) - \Gamma_j(x_j^k, x_j^k), y - w_j(x_j^k) \rangle \geq 0, \forall y \in C_j.$$

Let

$$y_j^k = x_j^k + \rho_j^{l(k,j)}(w_j(x_j^k) - x_j^k)) \text{ and } G_j(x_j^k) = F(y_j^k),$$

where $l(k,j)$ is the smallest integer which satisfies

$$\langle F(x_j^k + \rho_j^{l(k,j)}(w_j(x_j^k) - x_j^k)), x_j^k - w_j(x_j^k) \rangle \geq \alpha_j \| w_j(x_j^k) - x_j^k \|^2$$

**Step 6**. (cutting plane generation)

Set

$$A_j^{k+1} = \begin{bmatrix} A_j^k \\ G(x_j^k)^T \end{bmatrix}, \quad b_j^{k+1} = \begin{bmatrix} b_j^k \\ G(x_j^k)^T x_j^k \end{bmatrix},$$

$H_j^k = \{x \in \mathbb{R}^n; G(x_j^k)^T (x - x_j^k) = 0\}$ is the new cutting plane for $VI(F, C_j^k)$.

Increase $k$ by one GOTO Step 2.

By Theorem 1 of [2], similar to the proof of Theorem 2 we have

**Theorem 4** Let $F: X \to \mathbb{R}^n$ be Lipschitz continuous, i.e., there exists a constant $L > 0$ such that

$$(F(y) - F(x))^T (y - x) \leq L\|y - x\|, \forall x, y \in X$$

on a compact convex body $X$, and $X_D^*$ be nonempty. Then Algorithm 3 either stops with a solution of $VI(F, X)$ after a finite number of iterations, or there exists a subsequence of the infinite sequence $\{x_j^k\}$ that converges to a point $x^* \in X^*$. □

**Algorithm 4.** Step 1, Step 2, Step 4, Step 5, and Step 6 are same as those in Algorithm 3.

**Step 3**. (stop criterion)

Compute $g_X(x_j^k) = \max_{x \in X} F(x_j^k)^T (x_j^k - x)$,

if $g_X(x_j^k) < \varepsilon$, **then** STOP,

        **else** GOTO step 4.

By Theorem 4 we have

**Theorem 5** Let $F : X \to \mathsf{R}^n$ be Lipschitz continuous on a compact convex body $X$ and $X_D^*$ be nonempty. Then Algorithm 4 stops with a $\varepsilon$-solution of $VI(F, X)$ after a finite number of iterations.□

## 5. Generalized Analytic Center Cutting Plane Algorithms for Variational Inequalities with Unbounded Domains

This section presents an analytic center cutting plane algorithms for solving a strongly pseudomonotone variational inequality $VI[F, X]$ whose domain is an unbounded convex body (including n-dimensional Euclidean space of $\mathsf{R}^n$ itself). By use of Proposition 1 and 2, due to the pseudomonotonicity, $VI[F, X]$ has an unique solution $x^*$ over $X$. Let $\{C_j\}$ be a sequence of polytopes that satisfies

$$C_j \subseteq C_{j+1} \ (j = 1, \cdots), \text{ and } (\bigcup_{j=1}^{\infty} C_j)^c = X,$$

then $VI[F, C_j]$ have an unique solution $x_j^*$ over $C_j$ ($j = 1, 2, \ldots$). We can always assume that $C_j$ contains all boundary points of $X$ (if there are any). Since the solution $x^*$ of $VI[F, X]$ is a fixed point, $x^*$ lies in $C_j$ if $j$ is large enough (say $j > k$), therefore by Lemma 2 $x_j^* = x^*$ ($j > k$).

The following Algorithm 5 is proposed here to find $x^*$.

**Algorithm 5.**

**Step 1**. (initialization)
$k = 0, j = 1, A^k = A_j, b^k = b_j, C_j^k = \{x \in \mathsf{R}^n; A_j^k x \leq b_j^k\}$;

**Step 2**. (find an $\varepsilon$-solution for $\varepsilon = \frac{1}{2^j}$)

Find an approximate analytic center $x_j^k$ of $C_j^k$,

Compute $g_{C_j}(x_j^k) = \max_{x \in C_j} F(x_j^k)^T (x_j^k - x)$,

if $g_{C_j}(x_j^k) < \frac{1}{2^j}$, **then** increase $j$ by 1 RETURN TO Step1,

**else** GOTO Step 3;

**Step 3**. (cut generation)

Set $A_j^{k+1} = \begin{bmatrix} A_j^k \\ F(x_j^k)^T \end{bmatrix}$, $b_j^{k+1} = \begin{bmatrix} b_j^k \\ F(x_j^k)^T x_j^k \end{bmatrix}$,

$H_j^k = \{x \in \mathsf{R}^n; F(x_j^k)^T (x - x_j^k) = 0\}$ is the new cutting plane for $VI(F, C_j^k)$.

Increase $k$ by one GOTO Step 2.

**Theorem 6** Let $F: X \to \mathsf{R}^n$ be strongly monotone on $X$, then Algorithm 5 either stops with a solution of $VI(F, X)$ after a finite number of iterations, or there exists a subsequence of the infinite sequence $\{x_j^k\}$ that converges to a point $x^* \in X^*$.

**Proof**. $F$ is strongly monotone on $X$ implies that there exists a constant $N > 0$ such that [8]

$$\max_{y \in X} F(x)^T (x - y) \geq N \|x - y\|^2, \forall x \in X.$$

Let $x_j^*$ be the unique solution of $VI[F, C_j]$ over $C_j$ ($j = 1, 2, \ldots$). Suppose all boundary points of $X$ (if there are any) are in $C_j$ ($j = 1, 2, \ldots$). Then, if $j$ is large enough (say $j > k$), by Lemma 2 we have

$$x_j^* = x^* \ (j > k).$$

If Algorithm 5 does not stop after finite number of iterations, then exists an infinite sequence $\{x_j\} \subseteq X$ with $x_j \in C_j$ such that

$$g_{C_j}(x_j) < \frac{1}{2^j} \quad (j = 1, 2, \ldots).$$

Hence

$$N\|x_j - x^*\|^2 = N\|x_j - x_j^*\|^2 \leq g_{C_j}(x_j) = \max_{y \in C_j} F(x)^T(x_j - y) < \frac{1}{2^j} (j > k)$$

which implies that $\{x_j\}$ is a bounded sequence. Therefore, $\exists$subsequence of $\{x_j\}$, which is convergent to $x^{**}$ in $X$. Similar to the proof of Theorem 2, $x^{**}$ is a solution for $VI[F, X]$, and so $x^{**} = x^*$. $\square$

We notice that, in the proof of Theorem 6, the key is that $\{x_j\}$ is a bounded subsequence. Therefore, similarly have the following Theorem 7.

**Threorem 7** Let $F: X \to \mathbb{R}^n$ be strongly pseudomonotone on $X$, then Algorithm 5 either stops with a solution of $VI(F, X)$ after a finite number of iterations, or there exists a subsequence of the infinite sequence in $X$ that converges to a point $x^* \in X$. $\square$

Theorem 6 and 7 state that Algorithm 5 can always stop and output an approximate solution after finite number of iterations.

## 6. Conclusion Remark

This paper works with variational inequalities whose feasible sets are bounded or unbounded convex bodies. We present some analytic center cutting plane algorithms that extend the algorithms proposed in [2, 3, 26]. We should mention that our approach can be used to extend many interior methods which are associated with polyhedral feasible regions, e.g., the algorithms given by [5, 6]. We can also extend some other algorithms for variational inequalities over polyhedral feasible sets [7, 27, 28].